\documentclass[onecolumn]{elsart3p}

\usepackage{latexsym,amssymb,amsmath,amsfonts,url,bbm}
\usepackage[comma,sort&compress]{natbib}

\newtheorem{theorem}{Theorem}[section]
\newtheorem{lemma}{Lemma}[section]
\newtheorem{proposition}{Proposition}[section]
\newtheorem{corollary}{Corollary}[section]

\def \qed {  \mbox{}\hfill $\Box$\vspace{1ex}}
\newcommand{\equald}{\stackrel{\mathfrak{D}}{=}}
\newcommand{\convd}{\stackrel{\mathfrak{D}}{\to}}
\newcommand{\eps}{\varepsilon}

\begin{document}

\bibliographystyle{abbrvnat}

\begin{frontmatter}
\title{Maxima of Dirichlet and triangular arrays of gamma variables}
\author[addr]{Arup Bose},
\ead{abose@isical.ac.in}
\ead[url]{http://www.isical.ac.in/\textasciitilde abose}
\author[addr]{Amites Dasgupta},
\ead{amites@isical.ac.in}
\author[addr]{Krishanu Maulik\corauthref{cor}}
\corauth[cor]{Corresponding author.}
\ead{krishanu@isical.ac.in}
\ead[url]{http://www.isical.ac.in/\textasciitilde krishanu}

\address[addr]{Statistics and Mathematics Unit, Indian Statistical Institute, 203 B.T.~Road, Kolkata 700108, India}

\begin{abstract}
Consider a rowwise independent triangular array of gamma random
variables with varying parameters. Under several different
conditions on the shape parameter, we show that the sequence of
row-maximums converges weakly after linear or power transformation.
Depending on the parameter combinations, we obtain both Gumbel and
non-Gumbel limits.

The weak limits for maximum of the coordinates of certain Dirichlet
vectors of increasing dimension are also obtained using the gamma representation.
\end{abstract}

\end{frontmatter}

\begin{section}{Introduction} \label{sec: intro}
Suppose  $\{Y_n\}$ is a sequence of i.i.d.\ random variables and
${M}^*_n=\max\{Y_{1}, \ldots, Y_{n}\}$. Necessary and sufficient
conditions for the weak convergence of ${M}^*_n$ under linear
normalisation are well known. See for example,
\cite{fisher:tippett:1928, gnedenko:1943, dehaan:1970}. In
particular, let  $Y_n$ be i.i.d.\ standard normal variables and let
$G$ denote the Gumbel distribution
$$G(x)=\exp(-e^{-x}).$$
Then \citep[cf.][Theorem 1.5.3]{leadbetter:lindgren:rootzen:1983},
$$\lim_{n\to \infty} P [{M}^*_n \leq c_n x + d_n]
=G(x),$$ where
\begin{equation} \label{eq: alpha}
c_n=\frac 1{\sqrt{2 \log n}} \qquad \text{and} \qquad d_n =
\sqrt{2\log n} - \frac{\log \log n + \log (4\pi)}{2 \sqrt{2 \log
n}}.
\end{equation}

Now let $(Y_{1n}, \dots, Y_{nn})$ be a triangular sequence of random
variables and let $M_n =\max\{Y_{1n}, \dots, Y_{nn}\}$. The question
of convergence of $M_n$ has been addressed under a variety of
conditions.

For example, let $Y_{in}$ be i.i.d.\ with
$Y_{1n}=\big(\sum_{1 \leq j \leq {\alpha_n}}
U_j-\alpha_n\mu\big)/(\sigma \alpha_n^{1/2})$, where $U_j$ are
i.i.d.\ with mean $\mu$ and standard deviation $\sigma$; $\alpha_n$ is a sequence of integers going to $\infty$. Assuming
that $U_j$ has a finite moment generating function in an open
interval containing the origin and
\begin{equation} \label{eq: alpha rate}
\frac{\alpha_n^{(R+1)/(R+3)}}{\log n} \to \infty,
\end{equation}
for some integer $R \geq 0$,
\cite{anderson:coles:husler:1997} showed that
$$\lim_{n\to \infty} P [M_n \leq c_n x+d_n]
=G(x)$$ for $c_n$ as in~\eqref{eq: alpha} and some suitable
sequences $d_n$.

\cite{nadarajah:mitov:2002} considered the maximums of triangular
array of binomial, negative binomial and discrete uniform variables.
The case of binomial triangular array is discussed with increasing
number of trials $m_n$ and fixed probability of success, $p$.

\cite{bose:dasgupta:maulik:2007} considered the row-maximum of a
triangular array with dependent rows. More precisely, for
$n$-dimensional multinomial random variable with equally likely
cells, the maximum of the coordinates converges to Gumbel law if
number of trials increases fast enough.

We consider $\{Y_{in}\}$ to be a triangular sequence such that, for each
$n$, $Y_{in}$ are i.i.d.\ random variables having
Gamma~($\alpha_n$,~$1$) distribution. Also let $\boldsymbol{X}_n =
(X_{1n}, \ldots, X_{nn})$ be an $n$-dimensional Dirichlet
distribution with parameters $\alpha_n$,~$\ldots$,~$\alpha_n$,
$\beta_n$ supported on the $n$-dimensional simplex
$\{\boldsymbol{x}: 0\leq\sum_{i=1}^n x_i\leq 1\}$ and with density
$$\frac{\Gamma (n \alpha_n + \beta_n)}{\Gamma (\alpha_n)^n \Gamma
(\beta_n)} \left( \prod_{i=1}^n x_i \right)^{\alpha_n} \left( 1 -
\sum_{i=1}^n x_i \right)^{\beta_n}.$$ We investigate the problem of
existence of the weak limits of the maxima
$$M_n = \max \{Y_{1n}, Y_{2n}, \ldots, Y_{nn}\}, \qquad \widetilde M_n = \max \{X_{1n}, X_{2n}, \ldots, X_{nn}\}$$ under linear or power transformation.

In Section~\ref{sec: gamma}, we study the behavior of $M_n$. When
$n \alpha_n \to \infty$, the limit of centered and scaled $M_n$ is Gumbel, see Theorems~\ref{thm: log by alpha 0}--\ref{thm: eta lt 1}. In particular, if $\alpha_n / \log n \to 0$, we can take the scaling to be $1$. If $n \alpha_n$ has a positive, finite limit, in Theorem~\ref{thm: n alpha c}, we show that $M_n$ itself has a non-Gumbel limit. Under the assumption $n
\alpha_n \to 0$, the linear transformation of $M_n$ does not converge. However, in Theorem~\ref{thm: n alpha 0}, we show that a power transformation leads to uniform limit.

In Section~\ref{sec: dirichlet}, $\widetilde M_n$ is taken up. When
$n \alpha_n + \beta_n \to \infty$ and $n \alpha_n$ converges to a
positive limit, the limit of centered and scaled $\widetilde M_n$ is still Gumbel. When $n \alpha_n$ and $n \alpha_n + \beta_n$ both have finite, positive limits, $\widetilde M_n$ itself converges, but to a non-standard limit,
cf. Theorem~\ref{thm: dir n alpha c beta finite}. When $n \alpha_n \to 0$ and $n \alpha_n / \beta_n$ converges in $[0, \infty]$, in
Theorem~\ref{thm: dir n alpha zero}, we show that a power
transformation of scaled $\widetilde M_n$ converges to a mixture of
uniform distribution on $(0, 1)$ and a point mass at $1$.
\end{section}

\begin{section}{Maximum of triangular array of Gamma random
variables} \label{sec: gamma} The centering and scaling depends on
the nature of the sequence $\alpha_n$. The first case is similar to
Proposition 2 of \cite{anderson:coles:husler:1997}. Throughout the article, we use $\convd$ to denote convergence in distribution.

\begin{theorem} \label{thm: log by alpha 0}
Assume that $\alpha_n/\log n \to \infty$. Then
$$\sqrt{\frac{2 \log n}{\alpha_n}}(M_n - \alpha_n  - b_n
\sqrt{\alpha_n}) \convd G,$$ where $b_n$ is the unique
solution, in the region $b_n \sim \sqrt{2 \log n}$, of
\begin{equation} \label{def bn}
\log z + \frac{1}{2} \log (2 \pi) + z \sqrt{\alpha_n} - \alpha_n
\log \left(1 + \frac{z}{\sqrt{\alpha_n}}\right) = \log n.
\end{equation}
\end{theorem}

Observe that in this case,  $Y_{1n}$ can be considered to be the
``sum'' of $\alpha_n$ many i.i.d.\ random variables, each of which
is distributed as unit Exponential random variable. This set up is
similar to that of Proposition 2 of
\cite{anderson:coles:husler:1997} mentioned earlier  but we have the
added advantage that the random variables are gamma distributed. It
may also be noted that the condition $\alpha_n/\log n \to \infty$,
is the limiting form ($R=\infty$) of~\eqref{eq: alpha rate}. Almost
verbatim repetition of their argument yields the proof of
Theorem~\ref{thm: log by alpha 0}. We omit the details  but point
out that their Lemma~2 continues to hold if we replace the degree
$R$ polynomial in that lemma with the corresponding power series ($R=\infty$).
Using the moment generating function of the gamma distribution, the
$j$-th coefficient of the power series simplifies to $(-1)^{j+1} /
(j+2)$, for $j \geq 1$. This yields the defining equation for $b_n$
given in~\eqref{def bn} above.

Let the  centering and scaling required in general be $d_n$ and
$c_n$ respectively and let
\begin{equation} \label{def xn}
x_n = c_n x + d_n.
\end{equation}
Suppose $x_n$ is such that  $P[Y_{1n}> x_n]\to 0$. Then
\begin{equation} \label{def An}
-\log P[M_n \leq x_n]= -\log  P[Y_{1n} < x_n]^n\sim n P[Y_{1n}>
x_n].
\end{equation}

Motivated by the above, and noting that $Y_{1n}$ has the
Gamma$(\alpha_n, 1)$ distribution, define
$$A_n = \frac{n}{\Gamma(\alpha_n)} \int_{x_n}^\infty e^{-u}
u^{\alpha_n - 1} du = n P[Y_{1n} > x_n].$$

Integrating by parts, we immediately have
\begin{equation} \label{eq: abc}
A_n = B_n + (\alpha_n - 1) C_n,
\end{equation}
where
\begin{align}
B_n &= \frac{n}{\Gamma(\alpha_n)} e^{- x_n} x_n^{\alpha_n - 1}
\label{def Bn} \intertext{and} C_n &= \frac{n}{\Gamma(\alpha_n)}
\int_{x_n}^\infty e^{-u} u^{\alpha_n - 2} du \leq \frac{A_n}{x_n},
\label{def Cn}
\end{align}
which provides us with an upper bound for $A_n$:
\begin{equation} \label{ub An}
A_n \leq \frac1{1- \frac{|\alpha_n - 1|}{x_n}} B_n.
\end{equation}
For fixed $k$, if $\alpha_n > k$, we also obtain a lower bound via
integration by parts $k$ times repeatedly:
\begin{align}
A_n = & B_n \left[ 1 + \sum_{j=1}^{k-1} \prod_{i=1}^j
\frac{\alpha_n - i}{x_n} \right] + \prod_{i=1}^{k} (\alpha_n - i)
\frac{n}{\Gamma (\alpha_n)} \int_{x_n}^\infty e^{-u} u^{\alpha_n
-k -1} du \nonumber\\
\geq & B_n \left[ 1 + \sum_{j=1}^{k-1} \left(\frac{\alpha_n -
k}{x_n}\right)^j \right] = \frac{B_n}{1 - \frac{\alpha_n - k}{x_n}}
\left[1 - \left(\frac{\alpha_n - k}{x_n}\right)^k\right] \geq
\frac{B_n}{1 - \frac{\alpha_n - k}{x_n}} \left[1 -
\left(\frac{\alpha_n}{x_n}\right)^k\right]. \label{lb An}
\end{align}

When $\alpha_n / \log n$ remains bounded away from both $0$ and
$\infty$, we have $\alpha_n \to \infty$. The proof of the
following Theorem requires careful use of both the
bounds~\eqref{ub An} and~\eqref{lb An}.

\begin{theorem} \label{thm: alpha by log fin}
Assume that $\alpha_n/\log n$ remains bounded away from $0$ and
$\infty$. Then
$$\left(1-\frac{\alpha_n}{\zeta_n}\right)(M_n - \zeta_n) + \log
\left(1-\frac{\alpha_n}{\zeta_n}\right) \convd G,$$ where
$\zeta_n / \alpha_n$ is the unique solution bigger than $1$ of
\begin{equation} \label{def zetan}
z = 1 + \frac{\log n}{\alpha_n} - \frac{\log \sqrt{2 \pi} +
\frac12 \log \alpha_n}{\alpha_n} + \left(1 -
\frac1{\alpha_n}\right) \log z.
\end{equation}
\end{theorem}
\textbf{Proof.} We start with the solution of~\eqref{def zetan}. Observe that $(\log \sqrt{2 \pi} + \frac12 \log \alpha_n) / \alpha_n \to 0$, since
$\alpha_n \to \infty$. Further, since $\log n / \alpha_n$ is bounded
away from $0$, for all large $n$, $1 + \log n / \alpha_n - (\log
\sqrt{2 \pi} + \frac12 \log \alpha_n) / \alpha_n$ is bounded away
from  $1$. Thus,~\eqref{def zetan} has unique solution
$\zeta_n/\alpha_n$ bigger than 1 for all large $n$ and hence
$\zeta_n \to \infty$.

We have, $\zeta_n / \alpha_n > 1$, for all large $n$. If for some
subsequence $\{ n_k \}$, $\zeta_{n_{k}} / \alpha_{n_{k}} \to 1$,
then using~\eqref{def zetan}, we also have $\log n_k / \alpha_{n_k}
\to 0$, which is a contradiction. If for some subsequence $\{ n_k
\}$, $\zeta_{n_{k}} / \alpha_{n_{k}} \to \infty$, then
using~\eqref{def zetan}, we have
$$1 = \frac1{\zeta_{n_k} / \alpha_{n_k}} + \frac{\log n_k /
\alpha_{n_k}}{\zeta_{n_k} / \alpha_{n_k}} - \frac{\log \sqrt{2 \pi}
+ \frac12 \log \alpha_{n_k}}{(\zeta_{n_k} / \alpha_{n_k}) \cdot
\alpha_{n_k}} + \left( 1 - \frac1{\alpha_{n_k}} \right) \frac{\log
(\zeta_{n_k} / \alpha_{n_k})}{\zeta_{n_k} / \alpha_{n_k}}$$ and the
right side converges to $0$, which is again a contradiction. Thus,
$\zeta_n / \alpha_n$ is bounded away from both $1$ and $\infty$.

Corresponding to the choice of centering and scaling, we have
\begin{equation} \label{eq: xn n alpha 0}
x_n = \frac{x - \log \left(1 - \frac{\alpha_n}{\zeta_n} \right)}{1
- \frac{\alpha_n}{\zeta_n}} + \zeta_n.
\end{equation}
First observe that, since $\zeta_n / \alpha_n$ is bounded away
from both $1$ and $\infty$, we have $x_n \sim \zeta_n$. Hence,
\begin{equation} \label{eq: log xn n alpha 0}
\log x_n = \log \zeta_n + \frac{x - \log \left(1 -
\frac{\alpha_n}{\zeta_n} \right)}{\zeta_n \left( 1 -
\frac{\alpha_n}{\zeta_n} \right)} + O \left( \frac1{\zeta_n^2}
\right).
\end{equation}

Now, using Stirling's approximation,~\eqref{def Bn},~\eqref{eq: xn
n alpha 0},~\eqref{eq: log xn n alpha 0} and the fact $\alpha_n -
1 = O(\zeta_n) = o (\zeta_n^2)$, we have, after collecting terms,
\begin{align}
- \log B_n = & x \left[ \frac1{1 - \frac{\alpha_n}{\zeta_n}} -
\frac{\frac{\alpha_n - 1}{\zeta_n}}{1 - \frac{\alpha_n}{\zeta_n}}
\right] + \zeta_n - \frac{\log \left(1 - \frac{\alpha_n}{\zeta_n}
\right)}{1 - \frac{\alpha_n}{\zeta_n}} - \log n \nonumber\\
& - (\alpha_n -1) \log \zeta_n + \frac{\alpha_n - 1}{\zeta_n}
\frac{\log \left(1 - \frac{\alpha_n}{\zeta_n} \right)}{1 -
\frac{\alpha_n}{\zeta_n}} + \log \sqrt{2 \pi} \nonumber\\
& - \alpha_n + \left(\alpha_n - \frac12\right) \log \alpha_n +
o(1)
\nonumber\\
= & \alpha_n \left[\frac{\zeta_n}{\alpha_n} - 1 - \frac{\log
n}{\alpha_n} + \frac{\log \sqrt{2 \pi} + \frac12 \log
\alpha_n}{\alpha_n} - \left(1 - \frac1{\alpha_n}\right) \log
\frac{\zeta_n}{\alpha_n} \right] \nonumber\\
& + x \left[ 1 + \frac1{\zeta_n \left(1 - \frac{\alpha_n}{\zeta_n}
\right)} \right] - \frac{\log \left(1 - \frac{\alpha_n}{\zeta_n}
\right)}{1 - \frac{\alpha_n}{\zeta_n}} \left[ 1 -
\frac{\alpha_n}{\zeta_n} + \frac1{\zeta_n} \right] + o(1)
\label{eq: third expn}\\
= & x (1 + o(1)) - \log \left(1 - \frac{\alpha_n}{\zeta_n} \right)
+ o(1), \nonumber
\end{align}
where the first term of~\eqref{eq: third expn} vanishes since
$\zeta_n / \alpha_n$ satisfies~\eqref{def zetan}. Thus, we have
\begin{equation} \label{eq: Bn n alpha 0}
{\left(1 - \frac{\alpha_n}{\zeta_n} \right)}^{-1} B_n \sim e^{-x}.
\end{equation}
Hence, using the upper bound~\eqref{ub An} and the facts
$\alpha_n \to \infty$ and $x_n \sim \zeta_n$, we have
\begin{equation} \label{eq: upper bd}
\limsup A_n \leq \lim \frac{B_n}{1 - \frac{\alpha_n - 1}{x_n}} =
\lim \frac{B_n}{1 - \frac{\alpha_n}{\zeta_n}} = e^{-x}.
\end{equation}
Also, since $\alpha_n / \zeta_n$ is bounded away from $1$ and $x_n
\sim \zeta_n$, given any $\eps > 0$, we can fix a positive integer
$K$, such that $(\alpha_n / x_n)^K < \eps$, for all large $n$.
Hence, using the lower bound~\eqref{lb An}, since for all large $n$,
$\alpha_n > K$ and $(\alpha_n / x_n)^K < \eps$ hold, we have
$$\liminf A_n \geq \liminf \frac{B_n}{1 - \frac{\alpha_n - K}{x_n}} (1 -
\eps)= (1 - \eps) e^{-x}.$$ Since $\eps > 0$ is arbitrary, we get
$$\liminf A_n \geq e^{-x}.$$ Combining with~\eqref{eq: upper bd},
we get $A_n \to e^{-x}$, which completes the proof. \qed

When $n \alpha_n \to \infty$, but $\alpha_n = o(\log n)$, it turns
out that $c_n = 1$ and the limiting distribution is  $G$. However,
the choices of $d_n$ vary according to the specific limiting
behavior of $\alpha_n$. In general, we have the following lemma,
which is used repeatedly in the subsequent developments.
\begin{lemma} \label{lem: unified}
Suppose $c_n$ and  $d_n$ are such that for all $x \in \mathbb R$,
\begin{align}
B_n &\to - \log F(x), \label{eq: aim1} \intertext{and} (\alpha_n -
1) &= o (x_n). \label{eq: aim2}
\end{align}
Then $A_n\to -\log F(x)$ and hence
$$\frac{M_n - d_n}{c_n} \convd F.$$
\end{lemma}
\textbf{Proof.}
From the upper bound~\eqref{ub An}, we have
$$\limsup A_n \leq \lim_{n \to \infty} \frac1{1 - \frac{|\alpha_n
- 1|}{x_n}} \lim_{n \to \infty} B_n = -\log F(x).$$ Thus $A_n$ is
bounded. Also, since $(\alpha_n - 1)/x_n \to 0$, we have,
using~\eqref{def Cn},
$$|\alpha_n - 1| C_n \leq \frac{|\alpha_n - 1|}{x_n} A_n \to 0.$$

Then,~\eqref{eq: abc} gives us $\lim A_n = \lim B_n = -\log F(x)$.
\qed

As an illustration, suppose $\alpha_n=\alpha$ for all $n$.  It is
well-known that, in this case, the limiting distribution is Gumbel
with the centering, $d_n = \log n + (\alpha - 1) \log \log n - \log
\Gamma(\alpha)$, $c_n=1$  and $M_n - d_n \convd G$. See for
example, \citet[pp. 72-73]{resnick:1987}. This follows from the
above Lemma since $x_n = x + d_n \sim \log n$ and
$$- \log B_n =  - \log n + \log \Gamma(\alpha) + x_n -
(\alpha - 1) \log x_n = x + (\alpha - 1) \log \frac{\log n}{x_n}
\rightarrow  x.$$

We begin with the case where $\alpha_n \to \infty$, but $\alpha_n
= o (\log n)$.
\begin{theorem} \label{thm: alpha by log
0} Assume that $\alpha_n \to \infty$, such that $\alpha_n = o(\log
n)$. Then
$$M_n - \log n - (\alpha_n - 1) \log \log n - \xi_n + \log \Gamma(\alpha_n)
\convd G,$$ where $\xi_n / \log n$ is the unique positive
solution of
\begin{equation} \label{def xin}
z = \frac{\alpha_n - 1}{\log n} \log \left[ 1 - \frac{\alpha_n}{\log
n} \log \frac{\alpha_n}{\log n} + \frac{\alpha_n - \log \log n}{\log
n} + \frac{1}{2} \frac{\log \alpha_n}{\log n} + z \right].
\end{equation}
\end{theorem}
\textbf{Proof.}
First we consider the solution of~\eqref{def xin}. Define
\begin{align*}
\varepsilon_n &= \frac{\alpha_n - \log \log n}{\log n} +
\frac{1}{2} \frac{\log \alpha_n}{\log n} -  \frac{\alpha_n}{\log
n} \log \frac{\alpha_n}{\log n}\\
&= \frac{\alpha_n}{\log n} - \frac{\alpha_n - 1}{\log n} \log
\frac{\alpha_n}{\log n} - \frac12 \frac{\log \alpha_n}{\log n}\\
&\sim - \frac{\alpha_n}{\log n} \log \frac{\alpha_n}{\log n},
\end{align*}
since $\alpha_n \to \infty$, but $\alpha_n = o(\log n)$. Thus we
have $\eps_n > 0$ eventually, but $\varepsilon_n \to 0$. Also $m_n
:= \log n/(\alpha_n - 1) \to \infty$. With these
notations,~\eqref{def xin} becomes
$$e^{m_n z} = 1 + \eps_n + z.$$ Since $\eps_n > 0$ eventually,
there will be a unique positive solution $z_n$. For positive $z_n$,
we have $1 + \eps_n + z_n = e^{m_n z_n} > 1 + m_n z_n$, so that $z_n
< \eps_n / (m_n - 1) \to 0$. Hence,
\begin{equation} \label{eq: xi rate}
\xi_n = o(\log n).
\end{equation}

Using Stirling's formula, we write
\begin{align}
x_n =& x + d_n = x + \log n + (\alpha_n - 1) \log \log n + \xi_n -
\log \Gamma(\alpha_n) \label{def xn alpha by log 0}\\
=& x + \log n + (\alpha_n - 1) \log \log n + \xi_n - \log
\sqrt{2\pi} + \alpha_n \nonumber \\
& - (\alpha_n - 1/2) \log \alpha_n + o(1) \nonumber \\
=& \log n \Bigg[ 1 - \frac{\alpha_n}{\log n} \log
\frac{\alpha_n}{\log n} + \frac{\xi_n}{\log n} + \frac{\alpha_n
- \log \log n}{\log n} \nonumber \\
& + \frac{1}{2} \frac{\log \alpha_n}{\log n} + \frac{x - \log
\sqrt{2\pi}}{\log n} + o\left(\frac{1}{\log n}\right) \Bigg] \nonumber \\
=& R_n \log n \left[1 + \frac{x - \log \sqrt{2\pi}}{R_n\log n} + o
\left( \frac1{R_n \log n} \right) \right] \label{eq: Rn},
\end{align}
where
\begin{equation} \label{def Rn}
R_n = 1 - \frac{\alpha_n}{\log n} \log \frac{\alpha_n}{\log n} +
\frac{\alpha_n - \log \log n}{\log n} + \frac{1}{2} \frac{\log
\alpha_n}{\log n} + \frac{\xi_n}{\log n} \to 1,
\end{equation}
using~\eqref{eq: xi rate}. Hence
\begin{equation} \label{eq: d rate alpha by log 0}
x+d_n=x_n \sim d_n  \sim \log n
\end{equation}
and  $\alpha_n - 1 \sim \alpha_n = o(x_n)$, which gives
us~\eqref{eq: aim2}. Since $\xi_n/\log n$ is a solution
of~\eqref{def xin}, we have, using~\eqref{def Rn}, and
$\alpha_n/\log n\to 0$,
\begin{equation} \label{eq: xi R}
\xi_n = (\alpha_n - 1) \log R_n.
\end{equation}

Also, using~\eqref{def Bn},~\eqref{def xn alpha by log 0}
and~\eqref{eq: Rn}, we have
\begin{align*}
- \log B_n =& x + \xi_n - (\alpha_n - 1) \log \frac{x_n}{\log
n}\\
=& x + \xi_n - (\alpha_n - 1) \log R_n\\
&- (\alpha_n - 1)\log \left[1 + \frac{x - \log \sqrt{2\pi}}{R_n
\log n} + o\left(\frac{1}{R_n \log n}\right) \right]\\
=& x - (\alpha_n - 1)\log \left[1 + \frac{x - \log
\sqrt{2\pi}}{R_n \log n} + o\left(\frac{1}{R_n \log n}\right)
\right]\\
\sim& x - \frac{\alpha_n - 1}{\log n} (x - \log \sqrt{2\pi}) \to
x,
\end{align*}
which gives us~\eqref{eq: aim1} and completes the proof using
Lemma~\ref{lem: unified}.
\qed

Next we consider the cases where $\alpha_n$ is bounded above, but
$n \alpha_n \to \infty$. Here the centering $d_n$ depends on the
limiting behavior of $\log \alpha_n / \log n$. We separate out two
cases, depending on whether the ratio $\log \alpha_n / \log n$
converges to $0$, or is bounded away from $0$. We consider the
former case first, which includes the case $\alpha_n = \alpha$,
discussed earlier.

\begin{theorem} \label{thm: eta one}
Suppose $\alpha_n$ is bounded above, but $\log \alpha_n = o (\log
n)$. Then
$$M_n - \log n - (\alpha_n - 1) \log \log n + \log \Gamma
(\alpha_n) \convd G.$$
\end{theorem}
\textbf{Proof.}
In this case,
\begin{equation} \label{def xn eta one}
x_n = x+d_n=x + \log n + (\alpha_n - 1) \log \log n - \log \Gamma
(\alpha_n).
\end{equation}
If $\alpha_n \to 0$, then $\log \Gamma (\alpha_n) = - \log
\alpha_n + o(1) = o (\log n)$. Otherwise, $\alpha_n$ is bounded
away from both $0$ and $\infty$. Hence $\log \Gamma (\alpha_n)$ is
bounded and is $o(\log n)$. In either case, we have
\begin{equation} \label{eq: xn eta one}
x_n \sim \log n.
\end{equation}
Also note that
\begin{equation} \label{eq: d rate eta one}
d_n \sim x_n \sim \log n \sim \log (n \alpha_n).
\end{equation}
As $\alpha_n$ is bounded, we have $|\alpha_n - 1|= o(x_n)$, which
gives us~\eqref{eq: aim2}. Using~\eqref{def Bn} and~\eqref{def xn
eta one} and the fact $\alpha_n$ is bounded and~\eqref{eq: xn eta
one}, we have
$$- \log B_n = x - (\alpha_n - 1) \log \frac{x_n}{\log n} \to x.$$
This shows~\eqref{eq: aim1} and completes the proof of the theorem
using Lemma~\ref{lem: unified}.
\qed

Next we consider the case where $\alpha_n$ is bounded above and
$\log \alpha_n / \log n$ is bounded away from $0$.
\begin{theorem} \label{thm: eta lt 1}
Assume that $\alpha_n$ is bounded above, $n \alpha_n \to \infty$ and
$\log \alpha_n / \log n$ is bounded away from $0$. Then
$$M_n - \log (n \alpha_n) - (\alpha_n - 1) \log \log (n
\alpha_n) \convd G.$$
\end{theorem}
\textbf{Proof.}
From the given conditions,  we have $\log \alpha_n \to -\infty$
and hence $\alpha_n \to 0$. Here
\begin{equation} \label{def xn eta lt one}
x_n = x + d_n = x + \log (n \alpha_n) + (\alpha_n - 1) \log \log
(n \alpha_n).
\end{equation}
Since $n \alpha_n \to \infty$ and $\alpha_n \to 0$, we have
\begin{equation} \label{eq: xn eta lt one}
d_n \sim x_n \sim \log (n \alpha_n) \to \infty
\end{equation}
and thus, $\alpha_n - 1 = o (x_n)$, which gives us~\eqref{eq:
aim2}.

Also, using~\eqref{def Bn},~\eqref{def xn eta lt one},~\eqref{eq: xn
eta lt one} and the fact $\alpha_n \to 0$, we have
$$- \log B_n = x+\log \Gamma (\alpha_n+1) - (\alpha_n - 1) \log \frac{x_n}{\log (n \alpha_n)}
\to x,$$ Thus we have~\eqref{eq: aim1} and the proof is completed
using Lemma~\ref{lem: unified}.
\qed

Next we consider $\alpha_n$, which goes to $0$ at a faster rate.
We first look at the case $n \alpha_n \to \alpha \in (0, \infty)$.
In this case, the maximum $M_n$ itself converges to a
non-degenerate limiting distribution, which is parametrized by
$\alpha$. This distribution is not one of the three standard
classes of the extreme value distributions.
\begin{theorem} \label{thm: n alpha c}
Assume $n \alpha_n \to \alpha \in (0, \infty)$. Then, for all
$x>0$, we have
\begin{equation} \label{def lim n alpha c}
P [M_n \leq x] \to F_\alpha(x) := \exp \left( -\alpha
\int_x^\infty \frac{e^{-u}}u du \right), \quad x \geq 0.
\end{equation}
\end{theorem}
\textbf{Proof.}
The proof follows immediately from~\eqref{def An}. Using the
dominated convergence theorem, since  $P(Y_{1n} > x) \to 0$, we have
for all $x
> 0$,
$$- \log P [ M_n \leq x ] \sim n P [ Y_{1n} > x ] \sim n \alpha_n
\int_x^\infty e^{-u} u^{\alpha_n - 1} du \to \alpha \int_x^\infty
\frac{e^{-u}}u du.$$
\qed

When $n \alpha_n \to 0$, there does not exist any non-degenerate
limit distribution under linear transformations. However, a power
transformation gives Uniform $(0, 1)$ as the limiting
distribution. The idea behind the power scaling is contained in
the following lemma. This is used later in Section~\ref{sec:
dirichlet} as well.
\begin{lemma} \label{lem: gamma}
Let $V_n$ be Gamma $(\delta_n, 1)$ random variables, where $\delta_n
\to 0$. Then $V_n^{\delta_n} \convd U$, where $U$ is a Uniform
$(0, 1)$ random variable. Also, for all $k
> 0$, we have, $E[V_n^k] \to 1/(1+k)$.
\end{lemma}
\textbf{Proof.}
Observe that for any $k > 0$, we have $E \left[ V_n^{k \delta_n}
\right] = \Gamma (\delta_n (1+k)) / \Gamma (\delta_n) \sim 1/ (1+k)
= E [U^k]$. Then the result follows easily.
\qed
Thus, $Y_n^{\alpha_n}$ is approximately distributed as Uniform $(0,
1)$. Since the $n$-th power of the maximum of $n$ i.i.d.\ Uniform
$(0, 1)$ random variable is again Uniform $(0, 1)$, we expect
$M_n^{n \alpha_n}$ to converge to Uniform $(0, 1)$ distribution.

\begin{theorem} \label{thm: n alpha 0}
Assume that $n \alpha_n \to 0$. Then, for all $0<x<1$,
$$P [M_n^{n \alpha_n} \leq x] \to x.$$
\end{theorem}
\textbf{Proof.}
For any $0<x<1$, we have
\begin{align}
P [M_n^{n \alpha_n} \leq x] = & \left\{ \frac1{\Gamma (\alpha_n)}
\int_0^{x^{1/(n \alpha_n)}} e^{-u} u^{\alpha_n - 1} du \right\}^n
\nonumber\\
= & \left\{ \frac1{\alpha_n \Gamma (\alpha_n)} \int_0^{x^{1/n}}
e^{-u^{1/\alpha_n}} du \right\}^n \nonumber\\
= & \left\{ \frac{x^{1/n}}{\Gamma (\alpha_n + 1)} - \frac1{\Gamma
(\alpha_n + 1)} \int_0^{x^{1/n}} (1 - e^{-u^{1/\alpha_n}}) du
\right\}^n \nonumber\\
= & \frac{x}{\{\Gamma (\alpha_n + 1)\}^n} \left\{ 1 - \frac1n
\int_0^{x^{1/n}} \frac{n (1 - e^{-u^{1/\alpha_n}})}{x^{1/n}} du
\right\}^n. \label{eq: aim n alpha 0}
\end{align}

Since $\Gamma(x)$ is continuously differentiable on $(0, \infty)$,
we have
\begin{equation} \label{eq: gamma n alpha 0}
\{\Gamma (\alpha_n + 1)\}^n = (1 + O (\alpha_n))^n = 1 + O (n
\alpha_n) \to 1.
\end{equation}
We also have,
$$\int_0^{x^{1/n}} \frac{n (1 - e^{-u^{1/\alpha_n}})}{x^{1/n}} du
\leq \frac{n}{x^{1/n}} \int_0^{x^{1/n}} u^{1/\alpha_n} du = n
\alpha_n \frac{x^{1/(n \alpha_n)}}{\alpha_n + 1} \to 0,$$ since
$x^{1/(n \alpha_n)} \to 0$, for all $0 < x < 1$. Thus, we have,
$$\left\{ 1 - \frac1n \int_0^{x^{1/n}} \frac{n (1 -
e^{-u^{1/\alpha_n}})}{x^{1/n}} du \right\}^n \to 1.$$ The conclusion
follows from ~\eqref{eq: aim n alpha 0} and~\eqref{eq: gamma n alpha
0}.
\qed
\end{section}

\begin{section}{Maximum of coordinates of exchangeable Dirichlet
vectors with increasing dimension} \label{sec: dirichlet} We now extend the results to the maximum of Dirichlet distributions. The discussion is closely
related to the Gamma representation of Dirichlet: Recall
$\boldsymbol{X}_n$ is an $n$-dimensional vector having Dirichlet
$(\alpha_n, \ldots, \alpha_n; \beta_n)$ distribution. Let $\{Y_{in}:
1\leq i\leq n\}$ be i.i.d.\ Gamma ($\alpha_n$, $1$) random variables
and $Z_n$ be another independent Gamma ($\beta_n$, $1$) random
variable defined on the same probability space. Then
\begin{equation} \label{gamma reprn}
\boldsymbol{X}_n \equald \left( \frac{Y_{1n}}{\sum_{i=1}^n Y_{in}
+ Z_n}, \ldots, \frac{Y_{nn}}{\sum_{i=1}^n Y_{in} + Z_n} \right).
\end{equation}

This allows us to obtain the limiting distribution corresponding to
each case of i.i.d.\ Gamma random variables. For further
calculations, it helps to define $T_n = (\sum_{i=1}^n Y_{in} +
Z_n)/(n\alpha_n + \beta_n)$. Further, if $M_n = \max \{Y_{1n}, Y_{2n}, \ldots, Y_{nn}\}$, as before, then $(n\alpha_n + \beta_n) \widetilde M_n
\equald M_n / T_n$. So for some centering $d_n$ and scaling $c_n$,
we shall have
$$\frac{(n\alpha_n + \beta_n) \widetilde M_n - d_n}{c_n} \equald
\frac{M_n/T_n - d_n}{c_n} = \frac1{T_n} \left( \frac{M_n - d_n}{c_n}
- \frac{d_n}{c_n} (T_n-1) \right).$$ Note that, if $n \alpha_n +
\beta_n \to \infty$, we have $T_n \stackrel{P}{\to} 1$. Thus, we
have the following result as a simple application of Slutsky's
theorem, which we use repeatedly to obtain the results in Dirichlet
case.
\begin{proposition} \label{prop: dir}
Assume $n \alpha_n + \beta_n \to \infty$. Further assume that
$(M_n - d_n)/c_n \convd F$. Then $[(n \alpha_n + \beta_n)
\widetilde M_n - d_n]/c_n \convd F$, whenever
\begin{equation} \label{eq: dir cond}
\frac{d_n}{c_n} (T_n - 1) \stackrel{\mathrm{P}}{\to} 1.
\end{equation}
holds.
\end{proposition}

We now obtain the results for $\widetilde M_n$ as a corollary to
Proposition~\ref{prop: dir} above, whenever $n \alpha_n \to \alpha
\in (0, \infty]$ and $n \alpha_n + \beta_n \to \infty$.
\begin{corollary} \label{cor: dir}
Assume that $n \alpha_n \to \alpha \in (0, \infty]$ and $n
\alpha_n + \beta_n \to \infty$. Then
$$\frac{(n \alpha_n + \beta_n) \widetilde M_n - d_n}{c_n} \convd G.$$
with the centering
\begin{equation} \label{def dn}
d_n =
\begin{cases}
\alpha_n - b_n \sqrt{\alpha_n}, &\text{if $\log n =
o(\alpha_n)$},\\
\zeta_n - \left( 1 - \frac{\alpha_n}{\zeta_n} \right) \log \left(
1 - \frac{\alpha_n}{\zeta_n} \right), &\text{if
$\frac{\alpha_n}{\log n}$ is bounded away from both $0$ and $\infty$},\\
\log n + (\alpha_n - 1) \log \log n - \log \Gamma (\alpha_n) +
\xi_n, &\text{if $\frac{\alpha_n}{\log n} \to 0$ and $\alpha_n
\to \infty$},\\
\log n + (\alpha_n - 1) \log \log n - \log \Gamma (\alpha_n),
&\text{if $\log \alpha_n = o(\log n)$ and $\alpha_n$ bounded},\\
\log (n \alpha_n) + (\alpha_n - 1) \log \log (n \alpha_n),
&\text{if $\alpha_n$ is bounded, $n \alpha_n \to \infty$ and
$\frac{\log \alpha_n}{\log n}$}\\
&\quad \text{is bounded away from $0$},\\
0, &\text{if $n \alpha_n \to \alpha \in (0, \infty)$ and $\beta_n
\to \infty$},
\end{cases}
\end{equation}
and the scaling
\begin{equation} \label{def cn}
c_n =
\begin{cases}
\sqrt{\frac{\alpha_n}{2 \log n}}, &\text{if $\log n = o
(\alpha_n)$},\\
\left( 1 - \frac{\alpha_n}{\zeta_n} \right)^{-1}, &\text{if
$\frac{\alpha_n}{\log n}$ bounded from both $0$ and $\infty$},\\
1, &\text{if $\frac{\alpha_n}{\log n} \to 0$ and $n \alpha_n \to
\infty$},\\
\frac{\beta_n}{n \alpha_n + \beta_n}, &\text{if $n \alpha_n \to
\alpha \in (0, \infty)$ and $\beta_n \to \infty$},
\end{cases}
\end{equation}
where, as in Theorem~\ref{thm: log by alpha 0}, $b_n$ is the unique
root of~\eqref{def bn} in the region $b_n \sim \sqrt{2 \log n}$, as
in Theorem~\ref{thm: alpha by log fin}, $\zeta_n / \alpha_n$ is the
root of~\eqref{def zetan} which is bigger than $1$ and, as in
Theorem~\ref{thm: alpha by log 0}, $\xi_n / \log n$ is the unique
positive root of~\eqref{def xin}.

When $n \alpha_n \to \alpha \in (0, \infty)$ and $\beta_n \to
\infty$, the statement simplifies to
\begin{equation} \label{last conv}
\beta_n \widetilde M_n \convd F_\alpha.
\end{equation}
\end{corollary}
\textbf{Proof.}
We already have $(M_n - d_n)/c_n$ converges weakly to the
appropriate limit, $G$ or $F_\alpha$ from the corresponding
theorems in Section~\ref{sec: gamma}. We only verify~\eqref{eq:
dir cond} one by one by showing
$$Var \left[ \frac{d_n}{c_n} (T_n - 1) \right] = \left( \frac{d_n}{c_n}
\right)^2 \frac1{n \alpha_n + \beta_n} \to 0.$$

First consider $\log n = o (\alpha_n)$. From Theorem~\ref{thm: log
by alpha 0}, we have $d_n \sim \alpha_n$. Hence $d_n / c_n \sim
\sqrt{2 \alpha_n \log n}$. Thus,
$$\left( \frac{d_n}{c_n} \right)^2 \frac1{n \alpha_n + \beta_n}
\sim \frac{2 \alpha_n \log n}{n \alpha_n + \beta_n} \leq \frac{2
\log n}{n} \to 0.$$

Next we consider the cases where $\alpha_n \to \infty$ and
$\alpha_n / \log n$ is bounded above. We consider the cases where
the ratio converges to $0$ or stays bounded away from $0$
together. In the latter case, we know from Theorem~\ref{thm: alpha
by log fin} that $\alpha_n / \zeta_n$ is bounded away from $1$,
and thus, in~\eqref{def dn}, we have $d_n \sim \zeta_n =
O(\alpha_n) = O(\log n)$. Also, in~\eqref{def cn}, $c_n$ stays
bounded away from $0$ and $1$. Hence $d_n / c_n = O (\log n)$. If
$\alpha_n = o (\log n)$, then, by~\eqref{eq: d rate alpha by log
0}, we have $d_n \sim \log n$. Thus, in either case, $d_n / c_n =
O (\log n)$. Hence,
$$\left( \frac{d_n}{c_n} \right)^2 \frac1{n \alpha_n + \beta_n}
\sim \frac{\log^2 n}{n \alpha_n + \beta_n} \leq \frac{\log^2 n}{n}
\to 0.$$

Next we consider the cases, where $\alpha_n$ is bounded, $n
\alpha_n \to \infty$ and $\log \alpha_n / \log n$ is bounded
above. In these cases, $c_n = 1$. Using~\eqref{eq: d rate eta one}
and~\eqref{eq: xn eta lt one}, we have $d_n \sim \log (n
\alpha_n)$. Then,
$$\left( \frac{d_n}{c_n} \right)^2 \frac1{n \alpha_n + \beta_n}
\sim \frac{\log^2 (n \alpha_n)} {n \alpha_n + \beta_n} \leq
\frac{\log^2 (n \alpha_n)}{n \alpha_n} \to 0.$$

The remaining case, where $n \alpha_n \to \alpha \in (0, \infty)$
and $\beta_n \to \infty$ is trivial since $d_n = 0$. The
simplified form~\eqref{last conv} follows from Slutsky's theorem,
since $n \alpha_n + \beta_n \sim \beta_n$ in this case.
\qed

Next we consider the case,  $(n \alpha_n + \beta_n)$ remains bounded
and hence Proposition~\ref{prop: dir} does not hold.
\begin{theorem} \label{thm: dir n alpha c beta finite}
Assume that $n \alpha_n \to \alpha \in (0, \infty)$ and $\beta_n
\to \beta \in [0, \infty)$. Then
$$\widetilde M_n \convd H,$$ where $H$ is a distribution
supported on $(0, \infty)$ with $k$-th moment given by
${\mu_k}/\gamma_k$ where $\gamma_k$ is the $k$th moment of the Gamma
$(\alpha+\beta, 1)$ distribution and $\mu_k$ is the $k$-th moment of
the distribution $F_\alpha$, given by
$$\mu_k = \alpha \int_0^\infty x^{k-1} \exp \left( -x -\alpha \int_x^\infty
\frac{e^{-u}}u du \right) dx.$$
\end{theorem}
\textbf{Proof.}
From the representation~\eqref{gamma reprn} $\widetilde M_n \equald
M_n/S_n$, where $S_n = \sum_{i=1}^n Y_{in} + Z_n$. Further,
$({Y_{1n}}/{S_n}, \ldots, {Y_{nn}}/{S_n} )$ is independent of $S_n$
and hence  $M_n/S_n$ is independent of $S_n$.

From Theorem~\ref{thm: n alpha c}, we have $M_n \convd
F_\alpha$ and $S_n$ converges weakly to Gamma ($\alpha + \beta$)
distribution. Further since $\{M_n/S_n\}$ is bounded, it is tight.
Hence for any subsequence $n_k$ there is a further subsequence
$n_{k(l)}$, such that $\{M_{n_{k(l)}}/S_{n_{k(l)}}, S_{n_{k(l)}}\}$
converges weakly to say  $(V, W)$ where $W$ has Gamma ($\alpha +
\beta$) distribution. Since $M_n/S_n$ is independent of $S_n$, $V$
and $W$ are independent. Hence $VW$ has distribution $F_\alpha$.
This implies that $E(VW)^k=\mu_k$, that is $E(V^k)=\mu_k/\gamma_k$.
Since $V$ has support  $[0, 1]$, the moments determine the
distribution and the proof of the theorem is complete.
\qed

Now suppose $n \alpha_n \to 0$. In  this case, as for Gamma, no
linear transformation of $\widetilde M_n$ will have a limiting
distribution, however a power transformation will converge. We use
Lemma~\ref{lem: gamma} to obtain the following theorem.
\begin{theorem} \label{thm: dir n alpha zero}
Assume that $n \alpha_n \to 0$ and $n \alpha_n / \beta_n \to
\lambda \in [0, \infty]$. Then,
$$\left( \sigma_n \widetilde M_n \right)^{n \alpha_n} \convd
U_\lambda$$ where $U_\lambda$ is the distribution of $B_{\lambda} U
+ (1-B_{\lambda})$, $U$ and $B_{\lambda}$ are independent,
$P(B_{\lambda}=1)=\frac1{1+\lambda}=1-P(B_{\lambda}=0)$, $U$ is
uniform $(0, 1)$, and
$$\sigma_n =
\begin{cases}
\beta_n, &\text{when $\beta_n \to \infty$,}\\
1, &\text{otherwise.}\\
\end{cases}$$
When $\lambda=\infty$, we interpret $B_\infty$ as the random
variable degenerate at $0$.
\end{theorem}
\textbf{Proof.}
We first consider the case $\lambda=0$. Define $S_n = \sum_{i=1}^n
Y_{in} + Z_n$ as before.

When $\beta_n \to \infty$, clearly $S_n / (n \alpha_n + \beta_n)
\stackrel{\mathrm{P}}{\to} 1$ and $n \alpha_n + \beta_n \sim
\beta_n$. Thus, $S_n / \beta_n \stackrel{\mathrm{P}}{\to} 1$, and
hence $(S_n / \beta_n)^{n \alpha_n} \stackrel{\mathrm{P}}{\to} 1$.

If $\beta_n$ is bounded away from $0$ and $\infty$, $S_n$, which has
Gamma $(n \alpha_n + \beta_n)$ distribution, is  a tight sequence on
$(0, \infty)$ and hence $S_n^{n \alpha_n} \stackrel{\mathrm{P}}{\to}
1$.

If $\beta_n \to 0$, by Lemma~\ref{lem: gamma}, $S_n^{n \alpha_n +
\beta_n}$ converges to Uniform $(0,1)$ distribution weakly.  Hence
$S_n^{n \alpha_n} \stackrel{\mathrm{P}}{\to} 1$.

So whenever $\lambda=0$,  $S_n^{n \alpha_n}
\stackrel{\mathrm{P}}{\to} 1$. Also, from Theorem~\ref{thm: n alpha
0}, we have $M_n \convd U$. Thus by Slutsky's theorem and the
gamma representation~\eqref{gamma reprn}, we have the required
result when $\lambda=0$.

When $\lambda \in (0, \infty]$, the proof is very similar to that of
Theorem~\ref{thm: dir n alpha c beta finite}. When $\lambda =
\infty$, we shall interpret $\lambda / (1 + \lambda) = 1$ and $1 /
(1 + \lambda) = 0$. From Theorem~\ref{thm: n alpha 0}, we have that
$M_n^{n \alpha_n} \to U$. Also, from Lemma~\ref{lem: gamma}, we have
that $S_n^{n \alpha_n + \beta_n} \convd U$, and hence $S_n^{n
\alpha_n} \convd U^{\lambda/(1+\lambda)}$. Then arguing as in
Theorem~\ref{thm: dir n alpha c beta finite}, $\{M_n^{n \alpha_n} /
S_n^{n \alpha_n}\}$ converges weakly and the $k$-th moment of the
limit is given by
$$\frac{E[U^k]}{E\left[ U^{\lambda k/(1+\lambda)} \right]} =
\frac{1+\lambda+k}{(1+k)(1+\lambda)} = \frac\lambda{1+\lambda} +
\frac{1}{1+\lambda} \frac1{1+k},$$ which is the $k$-th moment of
the required limiting distribution.
\qed
\end{section}

\textbf{Acknowledgement.} We thank the Referee for careful reading and constructive suggestions.


\begin{thebibliography}{20}
\providecommand{\natexlab}[1]{#1}
\providecommand{\url}[1]{\texttt{#1}} \expandafter\ifx\csname
urlstyle\endcsname\relax
  \providecommand{\doi}[1]{doi: #1}\else
  \providecommand{\doi}{doi: \begingroup \urlstyle{rm}\Url}\fi

\bibitem[Anderson et~al.(1997)Anderson, Coles, and
  H{\"u}sler]{anderson:coles:husler:1997}
C.~W. Anderson, S.~G. Coles, and J.~H{\"u}sler.
\newblock Maxima of {P}oisson-like variables and related triangular arrays.
\newblock \emph{Ann. Appl. Probab.}, 7\penalty0 (4):\penalty0 953--971, 1997.
\newblock ISSN 1050-5164.

\bibitem[Bose et~al.(2007)Bose, Dasgupta, and Maulik]{bose:dasgupta:maulik:2007}
Arup Bose, Amites Dasgupta and Krishanu Maulik.
\newblock Maxima of the cells of an equiprobable multinomial.
\newblock \emph{Electron. Comm. Probab.}, 12:\penalty0 93--105 (electronic), 2007.
\newblock ISSN 1083-589X.

\bibitem[Fisher and Tippett(1928)]{fisher:tippett:1928}
R.~A. Fisher and L.~H.~C. Tippett.
\newblock Limiting forms of the frequency distribution of the largest and
  smallest members of a sample.
\newblock \emph{Proc. Camb. Philos. Soc.}, 24:\penalty0 180--190, 1928.

\bibitem[Gnedenko(1943)]{gnedenko:1943}
B.~Gnedenko.
\newblock Sur la distribution limite du terme maximum d'une s\'erie
  al\'eatoire.
\newblock \emph{Ann. of Math. (2)}, 44:\penalty0 423--453, 1943.
\newblock ISSN 0003-486X.

\bibitem[de~Haan(1970)]{dehaan:1970}
L.~de~Haan.
\newblock \emph{On regular variation and its application to the weak
  convergence of sample extremes}, volume~32 of \emph{Mathematical Centre
  Tracts}.
\newblock Mathematisch Centrum, Amsterdam, 1970.

\bibitem[Leadbetter et~al.(1983)Leadbetter, Lindgren, and
  Rootz{\'e}n]{leadbetter:lindgren:rootzen:1983}
M.~R. Leadbetter, G.~Lindgren, and H.~Rootz{\'e}n.
\newblock \emph{Extremes and related properties of random sequences and
  processes}.
\newblock Springer Series in Statistics. Springer-Verlag, New York, 1983.
\newblock ISBN 0-387-90731-9.

\bibitem[Nadarajah and Mitov(2002)]{nadarajah:mitov:2002}
S.~Nadarajah and K.~Mitov.
\newblock Asymptotics of maxima of discrete random variables.
\newblock \emph{Extremes}, 5\penalty0 (3):\penalty0 287--294, 2002.
\newblock ISSN 1386-1999.

\bibitem[Resnick(1987)]{resnick:1987}
S.~I. Resnick.
\newblock \emph{Extreme values, regular variation, and point processes}.
\newblock Springer-Verlag, New York, 1987.
\newblock ISBN 0-387-96481-9.

\end{thebibliography}
\end{document}